# From Manipulation to Abstraction: The Impact of Flexible Decomposition on Numerical Competence in Primary School

Fabio Pasticci

*Department of Philosophy, Social Sciences and Education, University of Perugia, Italy — fabio.pasticci@unipg.it*

**Abstract**

This study examines a structured instructional approach to the decomposition and recomposition of large numbers in six primary school classes (three 4th-grade and three 5th-grade, $N = 120$), using a quasi-experimental pre-test/post-test design with a control group and a four-week delayed post-test. The 12-week intervention follows the Concrete–Pictorial–Abstract (CPA) progression. The experimental groups achieved mean gains of 34.0 points (4th grade) and 29.6 points (5th grade) out of 100, against 16.4 and 11.1 points in the control groups (all $p < 0.001$). The Time × Group interaction in the mixed ANOVA gave partial $\eta^2 = 0.931$. An ANCOVA with the pre-test as covariate estimated an adjusted difference of 18.25 points, 95% CI [17.59, 18.92], $F(1, 117) = 2{,}978.10$, $p < 0.001$, partial $\eta^2 = 0.962$. The advantage was present on all five components of the test. At the delayed post-test the experimental groups retained 93.7% of their gain, against 78.6% in the control groups. Internal consistency of the instrument at pre-test was acceptable (Cronbach's $\alpha = 0.735$). The effects are unusually large, and the limits of the design that bear on their interpretation are discussed.

**Keywords:** place value; number decomposition; CPA approach; primary school; manipulatives

## 1. Introduction

Understanding the decimal number system is one of the fundamental competencies in primary school mathematics. The ability to decompose and recompose numbers is a prerequisite for mental arithmetic, written calculation and estimation (Fuson, 1990; Ross, 1989). Many pupils, however, have persistent difficulties with the principle of place value and tend to treat digits as separate entities rather than as components of an integrated system (Kamii, 1986).

### 1.1 Theoretical framework

The study rests on three references. Fuson (1992) describes a progression from concrete to symbolic, in which decomposition occupies the intermediate phase where the child coordinates the understanding of units with that of the positional system. Bruner's (1966) Concrete–Pictorial–Abstract framework translates this progression into instructional practice: manipulatives (base-ten blocks, abacus), then graphic representations, then symbolic notation. Sweller's (1988) cognitive load theory suggests why decomposition relieves working memory: breaking a complex task into manageable components, including non-canonical forms such as 1,200 = 12 hundreds, supports numerical flexibility, which is associated with success in mental arithmetic (Verschaffel et al., 2007).

## 1.2 Research questions

RQ1: To what extent does a structured CPA intervention improve performance on a test of place value and decomposition of large numbers, compared with ordinary instruction?

RQ2: Is the advantage uniform across the five components of the test, or concentrated in some of them?

RQ3: Is the advantage still present four weeks after the end of the intervention?

The design compares the intervention as a whole with ordinary instruction. It does not allow the contribution of individual elements, such as a particular decomposition strategy or a particular material, to be isolated; this is taken up among the limitations.

**2. Literature review**

## 2.1 Place value: cognitive challenges

Kamii (1986) showed that many children, even at the end of primary school, conceive digits as independent entities rather than as indicators of relative value. Ross (1989) identified five levels of place-value understanding and noted that many pupils stop at intermediate levels. Thompson and Bramald (2002) related this superficial understanding to difficulties in mental arithmetic with multi-digit numbers.

## 2.2 Decomposition and instructional approaches

Fuson and Briars (1990) documented the effectiveness of base-ten blocks, stressing the need for explicit guidance to connect concrete and symbolic representations. Carbonneau et al. (2013), in a meta-analysis of 55 studies, found an average effect size of $d = 0.37$ in favour of instruction with manipulatives. McNeil and Jarvin (2007) argued that concrete materials alone are not sufficient: pupils need explicit help in connecting concrete experience with abstraction.

## 2.3 Research gap

Few studies have compared structured CPA approaches with ordinary instruction in Italian primary schools with a specific focus on the decomposition of large numbers. This study contributes evidence from six classes and 120 pupils, with the raw data openly available.

**3. Method**

## 3.1 Participants and design

The study involved 120 pupils from six classes of a primary school in Italy: three 4th-grade classes (4A, 4B, 4C) and three 5th-grade classes (5A, 5B, 5C). Within each class, pupils were assigned to an experimental and a control subgroup of equal size (10, 9 and 11 pupils per subgroup in the three classes of each grade), giving 30 pupils in each of the four cells formed by group and grade.

The design is quasi-experimental, with a pre-test, a post-test immediately after the intervention, and a delayed post-test four weeks later. Because every class contributes pupils to both conditions, class is crossed with condition: whatever a class shares (teacher, timetable, classroom climate) affects both groups. The assignment of pupils to subgroups was not documented as a randomised procedure, so baseline equivalence cannot be assumed and is examined empirically (Sections 4.1 and 4.5). Splitting classes also carries a risk of contamination between the two subgroups, which would tend to reduce rather than inflate the observed difference.

## 3.2 Assessment instrument

The test consists of five parts, with a maximum total score of 100 points: (1) positional value (20 points); (2) canonical and polynomial decomposition (25 points); (3) reading and writing numbers (20 points); (4) ordering and comparison (25 points); (5) application in a real-world context (10 points). The same structure was used at the three time points.

Internal consistency, computed on the five part scores at pre-test, was acceptable (Cronbach's $\alpha = 0.735$). No other psychometric property of the instrument has been examined, and the total score is used here as given.

## 3.3 Instructional intervention

The intervention lasted 12 weeks and followed three phases. Concrete phase (weeks 1–4): base-ten blocks, physical abacus, Cuisenaire rods. Pictorial phase (weeks 5–8): schematic drawings, place-value diagrams, paper abacus. Abstract phase (weeks 9–12): standard notation, polynomial forms, non-canonical flexible decompositions. The control subgroups followed the ordinary curriculum with a symbolic approach.

## 3.4 Ethical considerations

The intervention and the tests were carried out by the class teachers as part of ordinary classroom activity and assessment. The author received the scores only, already anonymised, with no means of linking a score to a pupil's name or to an individual test paper. The archived dataset contains individual-level scores identified by anonymous codes and no names of pupils, teachers or schools.

## 3.5 Statistical analysis

The analysis comprises: (a) descriptive statistics with pooled within-class standard deviations; (b) paired-samples *t*-tests on individual gains; (c) a 2 (Group) × 2 (Grade) × 2 (Time) mixed ANOVA with repeated measures on Time; (d) an ANCOVA of post-test scores with the pre-test as covariate, after checking the homogeneity-of-slopes assumption; (e) a comparison of gains on the five parts of the test; (f) descriptive indices of retention at the delayed post-test. The significance threshold is $\alpha = 0.05$ (two-tailed) and effect sizes are reported as partial eta squared ($\eta^2 p$). Pupils are treated as independent observations; the consequences of this simplification are discussed in Section 5.4. All values were computed from the unrounded individual scores in the archived dataset (Pasticci, 2026).

**4. Results**

## 4.1 Descriptive statistics

Each cell (for example, Experimental 4th grade) is made up of three class subgroups of different sizes (10, 9 and 11 pupils). Their standard deviations are therefore combined by weighting each variance by its degrees of freedom:

$$SD_{pooled} = \sqrt{\frac{\sum \left( n_j - 1 \right) s_j^2}{\sum \left( n_j - 1 \right)}}$$

where $n_j$ and $s_j$ are the size and the standard deviation of class subgroup $j$, and $\sum \left( n_j - 1 \right) = 9 + 8 + 10 = 27$.

**Table 1.** Descriptive statistics of total scores (0–100), with pooled within-class standard deviations.

| Group | *N* | Pre-test *M* (*SD*) | Post-test *M* (*SD*) | Delayed *M* (*SD*) | Gain *M* | Gain % |
|---|---|---|---|---|---|---|
| Experimental 4th | 30 | 51.00 (1.99) | 85.00 (1.99) | 82.50 (2.00) | +34.00 | +66.7% |
| Control 4th | 30 | 49.80 (2.54) | 66.23 (2.01) | 63.20 (2.03) | +16.43 | +33.0% |
| Experimental 5th | 30 | 59.43 (2.80) | 89.00 (1.74) | 87.50 (1.78) | +29.57 | +49.7% |
| Control 5th | 30 | 59.73 (2.17) | 70.87 (2.08) | 68.00 (2.12) | +11.13 | +18.6% |

*Note.* Gain *M* is the mean of the individual differences between post-test and pre-test. Gain % is the gain as a percentage of the pre-test mean.

At pre-test the two groups were close. In 5th grade the difference was 0.30 points in favour of the control group, $t(58) = -0.43$, $p = 0.667$. In 4th grade the experimental group started 1.20 points higher, a small difference that nevertheless reaches significance, $t(58) = 2.05$, $p = 0.044$. This is the reason for the baseline-adjusted analysis in Section 4.5. At post-test the experimental cells are well above the control cells. The pooled standard deviation of the experimental 5th-grade cell falls from 2.80 to 1.74, while that of the experimental 4th-grade cell is unchanged (1.99).

## 4.2 Paired-samples *t*-tests

The paired-samples *t*-test asks whether the mean individual gain, $d_i = post_i - pre_i$, differs from zero. With $M_{diff}$ and $SD_{diff}$ the mean and standard deviation of the gains, the statistic is $t = M_{diff} / SE$, with $SE = SD_{diff} / \sqrt{n}$ and 29 degrees of freedom in each cell. The 95% confidence interval is $M_{diff} \pm t^* \left( 29 \right) \times SE$, with $t^* \left( 29 \right) = 2.045$.

**Table 2.** Paired-samples *t*-tests, pre-test to post-test.

| Group | *N* | *M* diff | *SD* diff | *t* | *df* | 95% CI |
|---|---|---|---|---|---|---|
| Experimental 4th | 30 | 34.00 | 1.98 | 93.93 | 29 | [33.26, 34.74] |
| Control 4th | 30 | 16.43 | 2.79 | 32.29 | 29 | [15.39, 17.47] |
| Experimental 5th | 30 | 29.57 | 2.43 | 66.62 | 29 | [28.66, 30.47] |
| Control 5th | 30 | 11.13 | 2.67 | 22.80 | 29 | [10.13, 12.13] |

*Note.* All $p < 0.001$ (two-tailed).

All four cells improved significantly. This shows only that each gain differs from zero, not that the gains differ between groups. Descriptively, across grades the experimental group ($N$ = 60) gained 31.78 points, $t(59) = 78.51$, 95% CI [30.97, 32.59], and the control group 13.78 points, $t(59) = 28.06$, 95% CI [12.80, 14.77]. The formal comparison between groups is given by the analyses that follow.

## 4.3 Mixed ANOVA

In the 2 (Group) × 2 (Grade) × 2 (Time) mixed ANOVA, Group and Grade are between-subjects factors and Time (pre-test, post-test) is a within-subjects factor. Between-subjects effects are tested against the between-subjects error, $MS$ = 6.93 ($df$ = 116); effects involving Time are tested against the within-subjects error, $MS$ = 3.10 ($df$ = 116). Partial eta squared is $SS_{effect}/\left(SS_{effect}+SS_{error}\right)$; by Cohen's (1988) conventions, values above 0.01, 0.06 and 0.14 are small, medium and large.

**Table 3.** Mixed ANOVA.

| Source | *SS* | *df* | *MS* | *F* | *p* | $\eta^2$p |
|---|---|---|---|---|---|---|
| *Between subjects* | | | | | | |
| Group (Exp. vs. Ctrl) | 5,358.15 | 1 | 5,358.15 | 773.29 | < 0.001 | 0.870 |
| Grade (4th vs. 5th) | 2,733.75 | 1 | 2,733.75 | 394.54 | < 0.001 | 0.773 |
| Group × Grade | 17.07 | 1 | 17.07 | 2.46 | 0.119 | 0.021 |
| Error | 803.77 | 116 | 6.93 | | | |
| *Within subjects* | | | | | | |
| Time (Pre vs. Post) | 31,144.82 | 1 | 31,144.82 | 10,060.70 | < 0.001 | 0.989 |
| Time × Group | 4,860.00 | 1 | 4,860.00 | 1,569.92 | < 0.001 | 0.931 |
| Time × Grade | 355.27 | 1 | 355.27 | 114.76 | < 0.001 | 0.497 |
| Time × Group × Grade | 2.82 | 1 | 2.82 | 0.91 | 0.342 | 0.008 |
| Error | 359.10 | 116 | 3.10 | | | |

The effect of Time, $F(1, 116) = 10{,}060.70$, $\eta^2$p = 0.989, reflects the improvement of all pupils and mixes the effect of the intervention with whatever is common to both groups (ordinary teaching, practice with the test, maturation). The effect that bears on RQ1 is the Time × Group interaction, $F(1, 116) = 1{,}569.92$, $p < 0.001$, $\eta^2$p = 0.931: the experimental group improved significantly more than the control group. The Time × Group × Grade interaction is not significant, $F(1, 116) = 0.91$, $p = 0.342$. The data therefore give no evidence that the effect differs between the two grades; with 60 pupils per grade, this is not evidence that it is the same.

## 4.4 Mixed ANOVA: explicit calculations

With two time points, the between-subjects sums of squares can be obtained from each pupil's mean of pre-test and post-test, and the within-subjects ones from each pupil's gain. The grand mean of the individual means is 66.38; the group means are 71.11 (experimental) and 61.66 (control); the grade means are 63.01 (4th) and 69.76 (5th). The grand mean of the gains is 22.78; by group, 31.78 and 13.78; by grade, 25.22 and 20.35. Each marginal mean is based on 60 pupils.

- $SS_{Group} = 2 \times 60 \times \left[ (71.11 - 66.38)^2 + (61.66 - 66.38)^2 \right] = 5,358.15$
- $SS_{Grade} = 2 \times 60 \times \left[ (63.01 - 66.38)^2 + (69.76 - 66.38)^2 \right] = 2,733.75$
- $SS_{Time} = 120 \times 22.78^2/2 = 31,144.82$
- $SS_{Time \times Group} = (60/2) \times \left[ (31.78 - 22.78)^2 + (13.78 - 22.78)^2 \right] = 4,860.00$
- $SS_{Time \times Grade} = (60/2) \times \left[ (25.22 - 22.78)^2 + (20.35 - 22.78)^2 \right] = 355.27$

The values shown use unrounded means; recomputing them from the two-decimal means above gives the same results to within rounding.

## 4.5 Analysis of covariance

Because the groups were not formed by a documented random procedure, and the 4th-grade groups differed slightly at pre-test, the post-test comparison is repeated with the pre-test as covariate:

$$Y_{post} = b_0 + b_1\,Group + b_2\,Y_{pre} + error$$

where Group is coded 0 for control and 1 for experimental, $b_1$ is the difference between groups at equal pre-test score, and $b_2$ is the common slope. In the full sample $b_0 = 44.74$ and $b_2 = 0.43$.

The model assumes that the regression lines of the two groups are parallel. The Group × Pre-test interaction was not significant in the full sample, $F(1, 116) = 0.108$, $p = 0.743$, in 4th grade, $F(1, 56) = 1.771$, $p = 0.189$, or in 5th grade, $F(1, 56) = 0.067$, $p = 0.797$, so the assumption is retained.

**Table 4.** ANCOVA: adjusted difference, adjusted means and effect size.

| Analysis | *N* | $b_1$ (adj. diff.) | 95% CI | *M* adj. Exp | *M* adj. Ctrl | *F* (*df*) | *p* | $\eta^2 p$ |
|---|---|---|---|---|---|---|---|---|
| Full sample | 120 | 18.25 | [17.59, 18.92] | 86.90 | 68.65 | 2,978.10 (1, 117) | < 0.001 | 0.962 |
| 4th grade | 60 | 18.41 | [17.42, 19.39] | 84.82 | 66.41 | 1,399.23 (1, 57) | < 0.001 | 0.961 |
| 5th grade | 60 | 18.23 | [17.30, 19.16] | 89.05 | 70.82 | 1,533.51 (1, 57) | < 0.001 | 0.964 |

*Note.* Adjusted means are evaluated at the mean pre-test score of the sample analysed in each row: 54.99 (full sample), 50.40 (4th grade) and 59.58 (5th grade).

At equal pre-test score, the experimental group scored 18.25 points higher than the control group, and the estimate is nearly the same in the two grades. In a model that also includes Grade and the Group × Grade interaction, the interaction is not significant, $F(1, 115) = 0.06$, $p = 0.80$, and the adjusted difference averaged over grades is 18.31 points ($SE = 0.33$); the covariate slope in that model is 0.31, since the slope of 0.43 in the simpler model also reflects the difference between grades. The mixed ANOVA and the ANCOVA rest on different assumptions and lead to the same conclusion about RQ1.

## 4.6 Components of the test

Table 5 compares the gains of the two groups on each of the five parts of the test, across grades.

**Table 5.** Mean gains (pre-test to post-test) on the five parts of the test, by group ($N = 60$ per group).

| Part (maximum) | Gain Exp | Gain Ctrl | Difference | *t*(118) | *p* |
|---|---|---|---|---|---|
| 1. Positional value (20) | 6.05 | 2.15 | 3.90 | 12.92 | < 0.001 |
| 2. Decomposition (25) | 7.38 | 3.35 | 4.03 | 13.50 | < 0.001 |
| 3. Reading and writing numbers (20) | 6.57 | 2.45 | 4.12 | 17.60 | < 0.001 |
| 4. Ordering and comparison (25) | 9.03 | 4.50 | 4.53 | 16.76 | < 0.001 |
| 5. Real-world application (10) | 2.75 | 1.33 | 1.42 | 6.74 | < 0.001 |

*Note.* Independent-samples *t*-tests on individual gains; the five tests are not adjusted for multiplicity.

The experimental group gained more on every part. In relation to the maximum score, the difference is largest for reading and writing numbers (4.12 of 20 points) and positional value (3.90 of 20), and smallest for real-world application (1.42 of 10). The advantage is therefore general rather than confined to the parts closest to the content of the intervention (RQ2), although it is somewhat weaker on the part that requires transfer to applied contexts.

## 4.7 Retention at four weeks

**Table 6.** Delayed post-test, four weeks after the post-test.

| Group | Post-test *M* | Delayed *M* | Loss | Delayed / Post | Share of gain retained |
|---|---|---|---|---|---|
| Experimental 4th | 85.00 | 82.50 | 2.50 | 97.1% | 92.6% |
| Control 4th | 66.23 | 63.20 | 3.03 | 95.4% | 81.5% |
| Experimental 5th | 89.00 | 87.50 | 1.50 | 98.3% | 94.9% |
| Control 5th | 70.87 | 68.00 | 2.87 | 96.0% | 74.3% |

*Note.* Share of gain retained = (Delayed − Pre) / (Post − Pre) × 100.

Both groups lost a few points between the post-test and the delayed post-test. The ratio of the delayed mean to the post-test mean differs little between groups (97.7% in the experimental group and 95.7% in the control group, across grades), because it carries the pre-test level in both numerator and denominator. The share of the gain that is retained separates the groups more clearly: 93.7% in the experimental group and 78.6% in the control group. In points, the mean loss was 2.00 in the experimental group and 2.95 in the control group, $t(118) = 9.57$, $p < 0.001$. The advantage of the experimental group is thus still present after four weeks (RQ3).

### 5. Discussion and conclusions

## 5.1 Main results

On RQ1, pupils who followed the CPA sequence improved much more than pupils who followed ordinary instruction: about 18 points more on a 100-point test, on which the sample averaged about 55 points at pre-test. The result holds in both grades and after adjustment for the pre-test.

The effect sizes ($\eta^2p$ = 0.931 for the Time × Group interaction and 0.962 in the ANCOVA) are far above the values usually reported in educational research, including the average of $d$ = 0.37 found by Carbonneau et al. (2013) for instruction with manipulatives. Effects of this size are uncommon and call for caution. They may reflect features of this study that favour large and homogeneous differences: a 12-week intervention closely aligned with the content of the test, small subgroups of 9 to 11 pupils, a test built for the study, and very low variability within cells. They should not be taken as an estimate of what the approach would produce in other schools.

Several mechanisms proposed in the literature are consistent with the direction of the result: a gradual passage from concrete to symbolic representations (Bruner, 1966), an extended manipulative phase, explicit guidance in connecting representations (McNeil & Jarvin, 2007), and exposure to several forms of decomposition (Verschaffel et al., 2007). The present design does not test any of them separately.

## 5.2 Components of the test and retention

On RQ2, the advantage appears on all five parts of the test, including those not limited to decomposition, such as ordering and comparison. It is smallest on real-world application, which suggests that transfer to applied contexts may need specific attention. On RQ3, the experimental group kept most of its gain after four weeks and lost slightly less than the control group.

## 5.3 What the study does not show

The data consist of scores on the five parts of the test. They do not record which decomposition strategies pupils used or which errors they made, so the study cannot say which strategy is most effective, nor how specific errors, such as those involving zero as a placeholder, changed. Nor can it separate the contribution of manipulatives from that of the other elements of the intervention. These questions require item-level or observational data.

## 5.4 Limitations

Eight limitations should be kept in mind. (1) Pupils were not assigned to subgroups by a documented random procedure; the pre-test adjustment reduces but does not remove the possibility of selection effects. (2) The two subgroups of each class worked in the same school context, so contamination between them cannot be excluded. (3) Pupils were analysed as independent observations although they belong to classes; class is crossed with condition, which limits the problem, but with three classes per grade the results describe these six classes rather than classes in general. (4) Teachers knew which pupils were in the experimental subgroup, and expectancy or novelty effects are possible. (5) The test was built for the study and only its internal consistency at pre-test has been examined. (6) The sample comes from a single school. (7) The follow-up was limited to four weeks. (8) Transfer to other areas of mathematics was not assessed.

## 5.5 Conclusions and future directions

In these six classes, a structured CPA approach to the decomposition of large numbers was associated with substantially larger gains than ordinary instruction, on all components of the test and with good retention after four weeks. Given the size of the effects and the limits of the design, the result is best read as strong encouragement to test the approach under more demanding conditions: random assignment of whole classes, an independently validated instrument, item-level data on strategies and errors, and a longer follow-up.

## Data availability

The raw scores are openly archived on Zenodo under a CC BY 4.0 licence (Pasticci, 2026): https://doi.org/10.5281/zenodo.20330352.